\documentclass[12pt]{article}

\usepackage[english]{babel}
\usepackage{amssymb,amsmath}
\usepackage{hyperref}

\newcommand{\footremember}[2]{%
    \footnote{#2}
    \newcounter{#1}
    \setcounter{#1}{\value{footnote}}%
}
\newcommand{\footrecall}[1]{%
    \footnotemark[\value{#1}]%
}

\newtheorem{theorem}{Theorem }[section]

\newtheorem{proposition}[theorem]{Proposition}

\newcommand{\proof}{\noindent\textbf{Proof. }}
\newcommand{\qed}{\hspace*{\fill}$\Box$}

\providecommand{\keywords}[1]{\textbf{Keywords---} #1}

\newcommand{\cB}{\mathcal{B}}

\title{Non-existence of two types of partial difference sets}

\author{S. De Winter \footremember{MTU}{Michigan Technological University} \footnote{sgdewint@mtu.edu} \and E. Neubert \footrecall{MTU} \footnote{ejneuber@mtu.edu} \and Z. Wang \footrecall{MTU} \footnote{zeying@mtu.edu}}

\date{}

\begin{document}

\maketitle 

\abstract{In this note we prove the non-existence of two types of partial difference sets in Abelian groups of order $216$. This finalizes the classification of parameters for which a partial difference set of size at most $100$ exists in an Abelian group.}


\section{Introduction}
Let $G$ be a finite Abelian group of order $v$. Then $D$ is a $(v,k,\lambda,\mu)$-{\it partial difference set} (PDS) in $G$ if $D$ is a $k$-subset of $G$ with the property that the expressions $gh^{-1}$, $g,h\in D$, represent each non-identity element in $D$ exactly $\lambda$ times, and each non-identity element of $G$ not in $D$ exactly $\mu$ times. Further assume that $D^{(-1)}=D$ (where $D^{(s)}=\{g^s:g\in D\}$ ) and $e\notin D$, where $e$ is the identity of $G$, then $D$ is called a {\it regular} partial difference set.  A regular PDS is called {\it trivial} if $D\cup\{e\}$ or $G\setminus\{D\}$ is a subgroup of $G$.

In \cite{MA94} Ma presented  a table of parameters  for which the existence of a regular PDS with $k\leq 100$ in an Abelian group was known or could not be excluded. In particular the list contained 32 cases where (non)-existence was not known. In \cite{MA97} Ma excluded the existence of a PDS in 13 of these 32 cases. In \cite{Klin} and \cite{Koh} existence was proved in one of the remaining cases, and recently De Winter, Kamischke and Wang \cite{SDWEKZW} proved nonexistence for all but two of the remaining cases. These remaining cases were the possible existence of a $(216, 40, 4, 8)$-PDS and a $(216, 43, 10, 8)$-PDS in an Abelian group of order $216$. In this note we will prove nonexistence of such PDS, hence finalizing the classification of parameters for which a PDS with $k\leq100$ exists in an Abelian group.  The proof uses ideas developed in \cite{SDWEKZW}, but requires an additional argument based on weighing points and lines in a projective plane. 

\section{Preliminaries}

The following three results will be used in our proof. The first two are due to Ma \cite{MA94,MADiscrete}, the third is a recent local multiplier theorem from \cite{SDWEKZW}.

\begin{proposition}\label{Ma2}
No non-trivial PDS exists in
\begin{itemize}
\item an Abelian group $G$ with a cyclic Sylow-$p$-subgroup and $o(G)\neq p$;
\item an Abelian group $G$ with a Sylow-$p$-subgroup isomorphic to $\mathbb{Z}_{p^s}\times\mathbb{Z}_{p^t}$ where $s\neq t$.
\end{itemize}
\end{proposition}

\begin{proposition}\label{Ma1}
Let D be a nontrivial regular $(v,\,k,\,\lambda,\,\mu)$-PDS in an Abelian group G. Suppose $\Delta=(\lambda-\mu)^2+4(k-\mu)$ is a perfect square.  If $N$ is a subgroup of $G$ such that $\gcd(\left|N\right|, \left|G\right|/\left|N\right|)=1$ and $\left|G\right|/\left|N\right|$ is odd, then $D_1=D\cap N$ is a (not necessarily non-trivial)  regular $(v_1,\,k_1,\,\lambda_1,\,\mu_1)$-PDS with
 $$\left|D_1\right|=\frac{1}{2}\left[ \left|N\right| +\beta_1\pm \sqrt{(\left|N\right|+\beta_1)^2-(\Delta_1-\beta_1^2)(\left|N\right|-1)}  \right].$$
Here $\Delta_1=\pi^2$ with $\pi=\gcd(\left|N\right|,\sqrt{\Delta})$ and $\beta_1=\beta-2\theta\pi$ where $\beta=\lambda-\mu$ and $\theta$ is the integer satisfying $(2\theta-1)\pi\leq\beta<(2\theta+1)\pi$.
\end{proposition}

\begin{proposition}\label{lmt}\rm{[LMT]}
Let D be a regular $(v,k,\lambda,\mu)$-PDS in an Abelian group $G$. Furthermore assume $\Delta=(\lambda-\mu)^2+4(k-\mu)$ is a perfect square.  Then $g\in G$ belongs to $D$ if and only if $g^s\in D$ for all $s$ coprime with $o(g)$, the order of $g$.
\end{proposition}

\section{The Main Result}

\begin{theorem}
There does not exist a $(216, 40, 4, 8)$-PDS in an Abelian group.
\end{theorem}
\proof  Assume by way of contradiction that $ D$ is a  $(216, 40, 4, 8)$-PDS in an Abelian group $G$ of order 216. By Proposition \ref{Ma2}, we know that $G\cong\mathbb{Z}_2^3 \times\mathbb{Z}_3^3$. 

Let $g_1$, $g_2$, $\dots$, $g_{26}$ be all elements of order 3 in $G$, and let $\cB_{g_i}=\{ag_i\,|\,o(a)=1 \;\mbox{or}\; 2, \; ag_i\in  D\}$, and $B_{i}=|\cB_{g_i}|$, $i=1,2, \dots, 26$. That is, $B_{i}$ equals the number of elements in $ D$ whose fourth power equals $g_i$. 

Now observe that the LMT implies that raising elements to the fifth power provides a bijection between $\cB_{g_i}$ and $\cB_{g_i^2}$. Hence $|\cB_{g_i}|=|\cB_{g_i^2}|$.

\medskip

Let $N$ be the Sylow-$2$-subgroup of $G$. Using Proposition \ref{Ma1} we obtain that $|N\cap D|=0$ or $4$. First assume that $D$ contains no elements of order $2$. We see that $\Sigma_iB_{i}=40$ and $\Sigma_iB_{i}(B_{i}-1)=56$, where the latter equality follows as all $7$ elements of order $2$ are not in $ D$, and thus each have exactly $\mu=8$ difference representations.\

By relabeling the $g_i$ if necessary, we may assume that $C_{j}:=B_{2j-1}=B_{2j}$, for $j=1,2,\dots, 13$, and $C_{1}\geq C_{2}\geq \dots\geq C_{13}$. We now obtain 
\begin{equation}\label{40_Case1}
\Sigma_jC_{j}=20 \quad \mbox {and} \quad \Sigma_j C_{j}^2=48.
\end{equation}

 It is easy to check that the system of equations (\ref{40_Case1}) exactly has the following nonnegative integer solutions, listed as 13 tuples $(C_{1},C_{2},\hdots,C_{13})$:
 \medskip

$(5,3,2,1,1,1,1,1,1,1,1,1,1)$, \;$(5,2,2,2,2,1,1,1,1,1,1,1,0)$,\

 $(4,4,2,2,1,1,1,1,1,1,1,1,0)$, \,$(4,3,3,2,2,1,1,1,1,1,1,0,0)$, \

$(4,3,2,2,2,2,2,1,1,1,0,0,0)$, \;$(4,2,2,2,2,2,2,2,2,0,0,0,0)$,\

$(3,3,3,3,2,2,1,1,1,1,0,0,0)$, \;$(3,3,3,2,2,2,2,2,1,0,0,0,0)$.
\medskip

 Secondly assume that $D$ contains $4$ elements of order $2$. It follows that $\Sigma_i B_{i}+4=40$. By counting the number of ways elements of order 2 can be written as differences of elements of $ D$, we obtain that  $\Sigma_iB_{i}(B_{i}-1)+4\cdot 3=4 \cdot 4+3 \cdot 8$. Using similar labeling as above, we now obtain 
\begin{equation}\label{40_Case2}
\Sigma_jC_{j}=18 \quad \mbox {and} \quad \Sigma_j C_{j}^2=32.
\end{equation}

 It is easy to check  that the system of equations (\ref{40_Case2}) has the following nonnegative integer solutions:
\medskip

$(3,3,2,1,1,1,1,1,1,1,1,1,1)$, \; $(3,2,2,2,2,1,1,1,1,1,1,1,0)$, \

$(2,2,2,2,2,2,2,1,1,1,1,0,0).$

\medskip

Recall that $N$ is the unique subgroup isomorphic to $\mathbb{Z}_2^3$ in $G$. Let $P_1,\hdots,P_{13}$ be the $13$ subgroups of $G$ isomorphic to $\mathbb{Z}_3$, and let $L_1,\hdots,L_{13}$ be the $13$ subgroups of $G$ isomorphic to $\mathbb{Z}_3^2$. Now consider the incidence structure $\mathcal{P}$ with points the subgroups $P_i\times N$, $i=1,\hdots,13$,  of $G$,  with blocks the subgroups $L_i\times N$, $i=1,\hdots,13$, of $G$, and with containment as incidence. Then it is easily seen that $\mathcal{P}$ is a $2-(13,4,1)$ design, or equivalently, the unique projective plane of order $3$. We next assign a weight to each point of $\mathcal{P}$ in the following way: if point $p$ corresponds to  subgroup $P_i\times N$ then the weight of $p$ is $\frac{1}{2}|((P_i\times N)\setminus N)\cap D|$. In this way the weights of the $13$ points of $\mathcal{P}$ correspond to the $13$ values $C_{1},C_{2},\hdots, C_{13}$, that is, half of the number of elements of order $3$ or $6$ from $ D$ in the subgroup underlying the given point. Without loss of generality we may assume the labeling is such that  point $P_i\times N$ has weight $C_i$. The weight of a block will simply be the sum of the weights of the points in that block.
\medskip

We next count how many elements of order $3$ or $6$ from $D$ a specific subgroup of the form $L_i\times N$ can contain. Assume that $|(L_i\times N) \cap D|=m$. Let $ag$ and $bh$ be two distinct elements from $ D$, with $a^2=b^2=g^3=h^3=e$. Then $agh^{-1}b^{-1}$ belongs to $L_i\times N$ if and only if $gh^{-1}\in L_i$. It is easy to see that if $g\in L_i$ there are $m-1$ possibilities for $bh$ such that $gh^{-1}\in L_i$, whereas if $g\notin L_i$ there are $\frac{| D|-m-2}{2}$ possibilities for $bh$ such that $gh^{-1}\in L_i$.
\medskip

 Counting the number of differences of elements of $ D$ that are in $L_i\times N$ in two ways, we obtain 

\begin{equation}\label{key}
m(m-1)+(k-m)(\frac{k-m-2}{2})=\lambda m +\mu(71-m),
\end{equation}
where $(k,\lambda,\mu)=(40, 4, 8)$. This yields that $m=8$ or $16$.
\medskip

Now define $m':=\frac{1}{2}|((L_i\times N)\setminus N) \cap D|$. We obtain the following table:
\medskip

\begin{tabular}{|c|c|c|}
\hline\
Case 1: $(216, 40, 4, 8)$-PDS & $D$ contains $0$ elements of order $2$  & $m'=4$ or $8$\\
\hline\
Case 2: $(216, 40, 4, 8)$-PDS & $D$ contains $4$ elements of order $2$ & $m'=2$ or $6$\\
\hline
\end{tabular}

\medskip

We now note that the values $m'$ must be the weights of the blocks of $\mathcal{P}$, and that in both cases these weights are even. We first show that no value $C_{i}$ can be odd. Assume by way of contradiction that $C_{i}$ is odd for some $i$. Let the weight of the four blocks that contain $P_i\times N$ be $n_{1},\hdots, n_{4}$ respectively. Then 

$$\sum_{j=1}^{13}C_{j}=C_{i}+\sum_{t=1}^{4}(n_t-C_{i}),$$

which implies that $\sum_{j=1}^{13}C_{j}$ is odd, contradicting with the fact that $\sum_{j=1}^{13}C_{j}=20$ or $18$.
\medskip

This leaves us with only the possibility $(4,2,2,2,2,2,2,2,2,0,0,0,0)$ for \linebreak[4] $(C_{1}\hdots,C_{13})$ in case 1. In this case, by considering the four blocks through a point with weight $2$ it easily follows that it is not possible to distribute the thirteen given weights in such a way that every block has weight $4$ or $8$. This concludes the proof. \qed

\bigskip

\begin{theorem}
There does not exist a $(216, 43, 10, 8)$-PDS in an Abelian group.
\end{theorem}

\proof This case is dealt with in a very similar way. We will only provide a sketch of the proof. Assume by way of contradiction $D$ is a $(216, 43, 10, 8)$-PDS in an Abelian group $G$.

As before $G\cong\mathbb{Z}_2^3 \times\mathbb{Z}_3^3$, and using Proposition \ref{Ma1} we obtain that $D$ contains either $3$ or $7$ elements of order $2$.
 If $D$ contains $3$ elements of order $2$ we obtain 
\begin{equation}\label{43_Case1}
\Sigma_jC_{j}=20\quad \mbox{ and} \quad \Sigma_j C_{j}^2=48
\end{equation}
 which is the same as the system of equations in (\ref{40_Case1}), and hence has the same set of solutions.\

 If $D$ contains $7$ elements of order $2$ we obtain 
\begin{equation}\label{43_Case2}
\Sigma_jC_{j}=18\quad \mbox{ and}\quad \Sigma_j C_{j}^2=32
\end{equation} which is the same as the system of equations in (\ref{40_Case2}),  and thus  has the same set of solutions.

With similar notation as in the previous theorem, and using the same counting argument for $(k,\lambda,\mu)=(43,10,8)$, one obtains $m=11$ or $19$, and
\medskip

\begin{tabular}{|c|c|c|}
\hline\
Case 3: (216, 43, 10, 8) & $D$ contains 3 elements of order 2  & $m'=4$ or $8$\\
\hline\
Case 4: (216, 43, 10, 8) & $D$ contains 7 elements of order 2  & $m'=2$ or $6$\\
\hline
\end{tabular}
\medskip

As before the weights of all blocks of $\mathcal{P}$ must be even, and the proof can be finished in the same way as in the $(216, 40, 4, 8)$-PDS case.\qed

\bigskip

\section{Conclusions}

It is interesting to note that no regular PDS exists in all but one of the cases that were originally left open in Ma's table \cite{MA94}. The exception arising from a two-weight code in an elementary Abelian $2$-group. Furthermore almost all known PDS in Abelian groups are of only few types: (negative) Latin square type, reversible difference sets, PCP type, Paley type, and projective two-weight sets. Also, recently it was shown that in Abelian groups of order $4p^2$, $p$ an odd prime, every non-trivial PDS is either of PCP type or a sporadic example in an Abelian group of order $36$ \cite{SDWZW}. These observations raise the question as to whether new strong and more general non-existence results can be proved, and whether further classifications for PDS in Abelian groups are possible. It is important to note that the situation in non-Abelian groups is very different, and many more examples exist in those groups.

\end{document}